\documentclass[12pt]{amsart}
\usepackage{amssymb}
\usepackage{amsfonts}
\usepackage{color,soul}
\usepackage{appendix}

\newcommand{\nc}{\newcommand}
\nc{\thref}[1]{Theorem~\ref{theo:#1}}
\nc{\selabel}[1]{\label{sect:#1}}
\nc{\seref}[1]{Section~\ref{sect:#1}}
\nc{\lelabel}[1]{\label{lemm:#1}}
\nc{\leref}[1]{Lemma~\ref{lemm:#1}}
\nc{\prlabel}[1]{\label{prop:#1}}
\nc{\prref}[1]{Proposition~\ref{prop:#1}}
\nc{\colabel}[1]{\label{coro:#1}}
\nc{\coref}[1]{Corollary~\ref{coro:#1}}
\nc{\exlabel}[1]{\label{exam:#1}}
\nc{\exref}[1]{Example~\ref{exam:#1}}
\nc{\delabel}[1]{\label{defi:#1}}
\nc{\deref}[1]{Definition~\ref{defi:#1}}
\nc{\eqlabel}[1]{\label{equa:#1}}
\nc{\relabel}[1]{\label{rema:#1}}
\nc{\reref}[1]{Lemma~\ref{rema:#1}}
\providecommand{\operatorname}[1]{\mathrm{#1}\,}
\nc{\Hom}{\operatorname{Hom}} \nc{\Mor}{\operatorname{Mor}}
\nc{\Aut}{\operatorname{Aut}} \nc{\Ann}{\operatorname{Ann}}
\nc{\Ker}{\operatorname{Ker}} \nc{\Trace}{\operatorname{Trace}}
\nc{\Char}{\operatorname{Char}} \nc{\Mod}{\operatorname{Mod}}
\nc{\End}{\operatorname{End}} \nc{\Spec}{\operatorname{Spec}}
\nc{\Span}{\operatorname{Span}} \nc{\sgn}{\operatorname{sgn}}
\nc{\Id}{\operatorname{Id}} \nc{\Com}{\operatorname{Com}}
\nc{\rank}{\operatorname{rank}}
\nc{\Clausen}{\operatorname{Cl}}
\nc{\Li}{\operatorname{Li}}
\nc{\Ls}{\operatorname{Ls}}

\let\:=\colon

\newtheorem{de}{Definition}[section]
\newtheorem{lm}[de]{Lemma}
\newtheorem{pr}[de]{Proposition}
\newtheorem{co}[de]{Corollary}
\newtheorem{re}[de]{Remark}
\newtheorem{res}[de]{Remarks}
\newtheorem{te}[de]{Theorem}
\newtheorem{ex}[de]{Example}
\newtheorem{exs}[de]{Examples}

\def\bex{\begin{ex}}
\def\eex{\end{ex}}
\def\bexs{\begin{exs}}
\def\eexs{\end{exs}}
\def\bl{\begin{lm}}
\def\el{\end{lm}}
\def\bc{\begin{co}}
\def\ec{\end{co}}
\def\bt{\begin{te}}
\def\et{\end{te}}
\def\bpr{\begin{pr}}
\def\epr{\end{pr}}
\def\br{\begin{re}}
\def\er{\end{re}}
\def\brs{\begin{res}}
\def\ers{\end{res}}
\def\bd{\begin{de}}
\def\ed{\end{de}}
\def\be{\begin{equation}}
\def\ee{\end{equation}}
\def\bea{\begin{eqnarray*}}
\def\eea{\end{eqnarray*}}
\def\bp{\begin{proof}}
\def\ep{\end{proof}}

\def\qed{\hfill\Box}

\let\:=\colon

\begin{document}

\title[\hfil Zagier formula for MZV's involving Hoffman elements]
{Another look at Zagier's formula for multiple zeta values involving Hoffman elements}

\begin{abstract}
In this paper, we give an elementary account into Zagier's formula for multiple zeta values involving Hoffman elements. Our approach allows us to obtain  direct proof in a special case via rational zeta series involving the coefficient $\zeta(2n)$. This formula plays an important role in proving Hoffman's conjecture which asserts that every multiple zeta value of weight $k$ can be expressed as a $\mathbb{Q}$-linear combinations of multiple zeta values of the same weight involving $2$'s and $3$'s. Also, using a similar hypergeometric argument via rational zeta series, we produce a new Zagier-type formula for the multiple special Hurwitz zeta values.  
\end{abstract}

\author{Cezar Lupu}



\maketitle

\section{Introduction}
\bigskip

Multiple zeta values, usually abbreviated MZV's (sometimes they are called Euler-Zagier sums), are real numbers defined by the convergent series

$$\displaystyle\zeta(k_{1}, k_{2}, \ldots, k_{r})=\sum_{0<n_{1}<n_{2}<\ldots <n_{r}}\frac{1}{n_{1}^{k_{1}}n_{2}^{k_{2}}\ldots n_{r}^{k_{r}}},$$
\bigskip
where $k_{1}, k_{2}, \ldots, k_{r}$ are positive integers with $k_{r}>1$. 
\vspace{0.3cm}

Originally defined by Euler for $r=2$, the theory of multiple zeta values was developed at the beginning of 90's independently by Hoffman \cite{Hoffman} and Zagier \cite{ZagierMZV,Zagier}. For more details about the theory of multiple zeta functions and their special values we recommend the monographs \cite{Burgos-Fresan, Zhao}. Also, the multiple zeta values can be defined as nested sums,

$$\displaystyle\zeta(k_{1}, k_{2}, \ldots, k_{r})=\sum_{n_{r}=1}^{\infty}\frac{1}{n_{r}^{k_{r}}}\sum_{n_{r-1}=1}^{n_{r}-1}\frac{1}{n_{r-1}^{k_{r-1}}}\ldots\sum_{n_{2}=1}^{n_{3}-1}\frac{1}{n_{2}^{k_{2}}}\sum_{n_{1}=1}^{n_{2}-1}\frac{1}{n_{1}^{k_{1}}}.$$
\vspace{0.3cm}

We call the above series a multiple zeta of depth $r$ and weight $k$, where $k=k_{1}+k_{2}+\ldots+ k_{r}$. Obviously, $0<r<k$ and there are $\binom{k-2}{r-1}$ multiple zeta values of given weight $k$ and depth $r$.
\vspace{0.3cm}

Although they look simple, it seems that these numbers have connections with Galois representation theory with applications in computing Feynman integrals in quantum field theory. Moreover, it turns out that MZV's satisfy $\mathbb{Q}$-linear combinations of powers of $\pi$ and odd values of the Riemann zeta function. In this regard, there are many conjectures concerning these values. One of them, which was conjectured by Hoffman \cite{Hoffman1} in 1997, and solved by Brown \cite{Brown} in 2012, asserts that that every multiple zeta value of a given weight $k$ can be expressed as a $\mathbb{Q}$-linear combination of MZV's of the same weight involving $2$'s and $3$'s.

The proof of this conjecture which is of motivic flavor relied on the evaluation of a
certain family of MZVs involving Hoffman elements. In other words, Brown \cite{Brown} showed that the following family of MZV's,
 
$$\displaystyle H(r, s)=\zeta(\underbrace{2, 2, \ldots, 2}_{\text{$r$}}, 3, \underbrace{2, 2, \ldots, 2}_{\text{$s$}} )$$
can be expressed as a $\mathbb{Q}$-linear combination of products $\pi^{2m}\zeta(2n+1)$, with $m+n=r+s+1$.

Unfortunately, Brown could not give an explicit formula for $H(r, s)$. This gap was filled by Zagier \cite{Zagier1} who provided the exact formula and its proof. 
\vspace{0.3cm}

In this paper, we give a different approach to Zagier's formula (for the special case $s=0$) using rational zeta series involving the values $\zeta(2n)$. The main idea is to look at the even powers in the Taylor series of $\arcsin$ and then to express these rational zeta series as a $\mathbb{Q}$-linear combination of powers of $\pi$ and $\zeta(2n+1)$. Moreover, a Zagier-type formula for multiple $t$-values is provided. In our approach, we use properties of special functions, e.g. Clausen function whose values are related to the values of the Riemann zeta function. More details are given in the next section.
\vspace{0.3cm}

This paper is organized as follows. In Section 2 we collect some basic tools about rational zeta series, Gauss hypergeometric function and Taylor series expansion of integer powers of $\arcsin$ which are used in the proof of the main result of the Section 2. Also, we give a proof of Zagier's formula for the special case $s=0$ using this rational zeta series approach. In Section 3, we give a Zagier type formula for the multiple $t$-values using the same approach as in Section 3. Finally, in the appendix, we provide some basic properties of the Clausen function which served as the main ingredient in the proof of Lemma 2.6. 
\vspace{0.3cm}

\textbf{Acknowledgements.} The author would like to thank Tom Hales, Camil Muscalu, Derek Orr for fruitful conversations which preceded this work, and to Don Zagier for inspiring conversations on the subject. Also, many thanks to my advisors Piotr Hajlasz and William C. Troy for encouragements through my PhD studies. This paper is part of author's PhD thesis at the University of Pittsburgh. 
\vspace{0.3cm}

\section{Zagier's formula for Hoffman elements}
\bigskip

The formula for $H(r, s)$ was provided by Zagier \cite{Zagier1} in the following
\vspace{0.3cm}

\bt
The following formula holds true:

$$\displaystyle H(r, s):=\zeta(\underbrace{2, 2, \ldots, 2}_{\text{$r$}}, 3, \underbrace{2, 2, \ldots, 2}_{\text{$s$}} )=2\sum_{k=1}^{r+s+1}(-1)^kc_{r, s}\zeta(2k+1)\zeta(\underbrace{2, 2, \ldots, 2}_{\text{$r+s+1-k$}}),$$

where $\displaystyle c_{r, s}^{k}=\binom{2k}{2r+2}-\left(1-\frac{1}{2^{2k}}\right)\binom{2k}{2s+1}$ and $r, s\geq 0$ are integers. 
\et
\vspace{0.3cm}

As we mentioned in the introduction, this formula served as the missing ingredient in Brown's proof of Hoffman's conjecture. Moreover, this played a key role in proving the zig-zag conjecture of Broadhurst \cite{Brown-Schnetz} by Brown and Schnetz. 

Also, it is worth to mention that Brown's theorem (Hoffman's conjecture) implies that if $\mathcal{Z}_{k}$ is the $\mathbb{Q}$-vector space spanned by all multiple zeta values of weight $k$, then

$$\displaystyle\dim\mathcal{Z}_{k}\leq d_{k},$$
where $d_{k}$ is the coefficient of $x^k$ in the powers series expansion of $\displaystyle\frac{1}{1-x^2-x^3}$.
\vspace{0.3cm}

The proof given by Zagier \cite{Zagier1} is elementary, but indirect. In other words, the proof of Theorem 2.1 is given by computing the associated generating functions of both sides in a closed form, and then showing they are entire functions of exponential growth that agree at sufficiently many points to force their equality. However, his proof was simplified by Li in \cite{Li}.

\bigskip
The special case when $s=0$ of Theorem 2.1 reads as following
\bigskip

\bt
We have

$$\displaystyle {H(r, 0):=\zeta(\underbrace{2, 2, \ldots, 2}_{\text{$r$}}, 3)=2\sum_{k=1}^{r+1}(-1)^kc_{r, 0}^k\zeta(2k+1)\zeta(\underbrace{2, 2, \ldots, 2}_{\text{$r+1-k$}})},$$

where $\displaystyle c_{r, 0}^k=\binom{2k}{2r+2}-\left(1-\frac{1}{2^{2k}}\right)2k$, and $r\geq 0$.

\et
\vspace{0.4cm}

Before we delve into the proof of the above theorem, let us recall the famous Euler's formula for $\zeta(2)$ which asserts that $$\displaystyle\zeta(2)=\sum_{n=1}^{\infty}\frac{1}{n^2}=\frac{\pi^2}{6}.$$

One of the many proofs of this result relies on the power series expansion of $\arcsin^2(x)$ near $x=0$,

\begin{equation}
    \arcsin^2(x)=\sum_{n=1}^{\infty}\frac{2^{2n-1}}{n^2\binom{2n}{n}}x^{2n}, |x|<1,
\end{equation}

and by using Wallis' integral formula after we substitute $x=\sin t$. More generally, if we integrate $(1)$ we obtain
\bigskip
$$\displaystyle\sum_{n=1}^{\infty}\frac{2^{2n}}{n^3\binom{2n}{n}}x^{2n}=4\int_0^y\frac{\arcsin^2x}{x}dx=2x^2\cdot {}_4F_{3}\left(1, 1, 1, 1; \frac{3}{2}, 2, 2; x^2\right),$$

where $\displaystyle {}_p F_{q}(a_1,a_2,\dots,a_p;b_1,b_2,\dots,b_q;z)=\sum_{n=0}^{\infty}\frac{(a_{1})_{n}\cdot\ldots (a_{p})_{n}}{(b_{1})_{n}\cdot\ldots (b_{q})_{n}}\cdot\frac{x^n}{n!}, |x|<1$,

and $(a)_{n}=a(a+1)\ldots(a+n-1)$ is the Pochhammer symbol.
\vspace{0.3cm}

For instance, when $y=\frac{1}{2}$ we have

$$\displaystyle\sum_{n=1}^{\infty}\frac{1}{n^3\binom{2n}{n}}=4\int_0^{\frac{1}{2}}\frac{\arcsin^2(x)}{x}dx=-2\int_0^{\frac{\pi}{3}}z\log\left(2\sin\frac{z}{2}\right)dz=$$

$$\displaystyle =-\frac{\zeta(3)}{3}-\frac{\pi\sqrt{3}}{72}\left(\Psi\left(\frac{1}{3}\right)-\Psi\left(\frac{2}{3}\right)\right),$$
\vspace{0.3cm}
where $\Psi(z)$ is the trigamma function. Also, for $y=1$ we have the special value

$$\displaystyle\int_0^1\frac{\arcsin^2x}{x}dx=\frac{\pi^2\log 2}{4}-\frac{7}{8}\zeta(3).$$
\bigskip

Moreover, using $(1)$, one can also derive Euler's 1775 famous formula \cite{Ayoub} for Apery's constant,

$$\displaystyle\zeta(3)=-2\pi^2\sum_{n=0}^{\infty}\frac{\zeta(2n)}{(2n+2)(2n+3)2^{2n}}.$$
\vspace{0.3cm}

Similarly, The Taylor series expansions for $\arcsin^4(x)$ and $\arcsin^6(x)$ near $x=0$ are given by
\vspace{0.3cm}

\begin{equation}
\arcsin^4x=\frac{3}{2}\sum_{n=1}^{\infty}\left(\sum_{m=1}^{n-1}\frac{1}{m^2}\right)\frac{1}{2^{2n}n^2\binom{2n}{n}}x^{2n},    
\end{equation}
and

\begin{equation}
\arcsin^6x=\frac{45}{4}\sum_{n=1}^{\infty}\left(\sum_{m=1}^{n-1}\frac{1}{m^2}\sum_{p=1}^{m-1}\frac{1}{p^2}\right)\frac{4^n}{\binom{2n}{n}n^2}x^{2n}. 
\end{equation}
\vspace{0.3cm}

As it has already been highlighted in \cite{Borwein-Chamberland, Chu-Zheng, Leshchiner, Sun} one can obtain the even integer powers for the Taylor series of arcsin by comparing the coefficients like powers of $\lambda$ in the formulas

$$\displaystyle\cos(\lambda\arcsin x)={}_{2} F_{1}\left(\frac{\lambda}{2}, -\frac{\lambda}{2}; \frac{1}{2}; x^2\right),$$
where ${}_2F_{1}$ is the Gauss hypergeometric function. This implies the following Taylor series expansion,

\begin{equation}
   \displaystyle\frac{\arcsin^{2r}(x)}{(2r)!}=\frac{1}{4^r}\sum_{n=1}^{\infty}\frac{4^n}{n^2\binom{2n}{n}}\cdot x^{2n}\cdot\sum_{n_{1}<n_{2}<\ldots n_{r-1}<n}\frac{1}{n_{1}^2n_{2}^2\ldots n_{r-1}^2}.
\end{equation}
\bigskip
\vspace{0.3cm}

Similarly, Chu and Zheng \cite{Chu-Zheng} gave explicit closed form for integer powers of $\arcsin$ in terms of elementary symmetric functions. In other words, let

$$\displaystyle e_{r}(x_{1}, x_{2}, \ldots, x_{r})=\prod_{i=1}^{r}(1+tx_{i})=\sum_{r} e_{r}(x_{1}, x_{2}, \ldots, x_{r})\cdot t^m.$$

By the same hypergeometric functions argument as above, they reproved $(4)$ in the following form,

\begin{equation}
    \displaystyle\sum_{n=1}^{\infty}\frac{(2x)^{2n}}{n^2\binom{2n}{n}}e_{r}\left(1, \frac{1}{2^2}, \ldots, \frac{1}{(n-1)^2}\right)=\frac{(2\arcsin x)^{2r+2}}{(2r+2)!}.
\end{equation}
\bigskip

\bigskip

Before we proceed with our main results, let us recall the much easier formula for the simplest of the Hoffman elements whose new proof is in the spirit of the ideas from this paper.

\bigskip

\begin{pr}[Hoffman-Zagier]
We have

$$\zeta(\underbrace{2, 2, \ldots, 2}_{\text{$r$}})=\frac{\pi^{2r}}{(2r+1)!},$$

\end{pr}
\vspace{0.3cm}

\textit{Proof.} Substituting $x=\sin t$ in formula $(4)$ we have

$$\displaystyle\frac{t^{2r}}{(2r)!}=\frac{1}{4^r}\sum_{n=1}^{\infty}\frac{4^n}{n^2\binom{2n}{n}}\sin^{2n}t\sum_{1\leq n_{1}<n_{2}<\ldots n_{r-1}<n}\frac{1}{n_{1}^2n_{2}^2\ldots n_{r-1}^2} $$
\vspace{0.3cm}

Integrating from $0$ to $\frac{\pi}{2}$, and using Wallis' formula, $\displaystyle\int_0^{\frac{\pi}{2}}\sin^{2n}tdt=\frac{\pi\binom{2n}{n}}{2^{2n+1}}$, we have

$$\displaystyle\frac{\pi^{2r}}{(2r+1)!}=\sum_{n=1}^{\infty}\frac{1}{n^2}\sum_{1\leq n_{1}<n_{2}<\ldots n_{r-1}<n}\frac{1}{n_{1}^2n_{2}^2\ldots n_{r-1}^2}=\zeta(\underbrace{2, 2, \ldots, 2}_{\text{$r$}}) $$

and we are done.$\qed{}$
\vspace{0.3cm}

\textbf{Remark.} Using $(5)$, the above formula can be expressed in terms elementary symmetric polynomials as
\vspace{0.3cm}

$$\displaystyle\sum_{n=1}^{\infty}\frac{1}{n^2}e_{r}\left(1, \frac{1}{2^2}, \ldots, \frac{1}{(n-1)^2}\right)=\frac{\pi^{2r+2}}{(2r+3)!}.$$
\bigskip

\bigskip

In what follows we will derive Theorem 2.2 in a direct and more elementary way using rational zeta series involving the coefficient $\zeta(2n)$. We have
\bigskip

\bt
We have

$$\displaystyle H(r, 0):=\zeta(\underbrace{2, 2, \ldots, 2}_{\text{$r$}}, 3)=-4(2r+3)\sum_{n=0}^{\infty}\frac{\zeta(2n)}{(2n+2r+2)(2n+2r+3)2^{2n}}\zeta(\underbrace{2, 2, \ldots, 2}_{\text{$r+1$}}).$$

\et

\bigskip

\textit{Proof of Theorem 2.4}. Dividing by $x$ and integrating from $0$ to $\sin t$, we have
\bigskip

$$\displaystyle \int_0^{\sin t}\frac{\arcsin^{2r}(x)}{x}dx=\frac{(2r)!}{4^r}\sum_{n=1}^{\infty}\frac{4^n}{2n^3\binom{2n}{n}}\sin^{2n}t\sum_{n_{1}<n_{2}<\ldots n_{r-1}<n}\frac{1}{n_{1}^2n_{2}^2\ldots n_{r-1}^2}.$$
\bigskip

By the substitution $x=\sin u$ in the integral from the left hand side, we have 

$$\displaystyle\int_0^{\sin t}\frac{\arcsin^{2r}(x)}{x}dx=\int_0^tu^{2r}\cot udu=\frac{(2r)!}{4^r}\sum_{n=1}^{\infty}\frac{4^n}{2n^3\binom{2n}{n}}\sin^{2n}t\sum_{n_{1}<n_{2}<\ldots n_{r-1}<n}\frac{1}{n_{1}^2n_{2}^2\ldots n_{r-1}^2}. $$

On the other hand, using $\displaystyle u\cot u=-2\sum_{n=0}^{\infty}\frac{\zeta(2n)}{\pi^{2n}}u^{2n}$, $|u|<\pi$, we derive

$$\displaystyle\int_0^tu^{2r}\cot udu=-2\int_0^t\sum_{n=0}^{\infty}\frac{\zeta(2n)}{\pi^{2n}}u^{2n+2r-1}du=-2\sum_{n=0}^{\infty}\frac{\zeta(2n)}{\pi^{2n}}\cdot\frac{t^{2n+2r}}{(2n+2r)}.$$
\bigskip

Therefore, we obtain

$$\displaystyle -2\sum_{n=0}^{\infty}\frac{\zeta(2n)}{\pi^{2n}}\cdot\frac{t^{2n+2r}}{(2n+2r)}=\frac{(2r)!}{4^r}\sum_{n=1}^{\infty}\frac{4^n}{2n^3\binom{2n}{n}}\sin^{2n} t\sum_{n_{1}<n_{2}<\ldots n_{r-1}<n}\frac{1}{n_{1}^2n_{2}^2\ldots n_{r-1}^2}.$$
\bigskip

Finally, integrating from $0$ to $\frac{\pi}{2}$ and using Wallis' integral formula, $\int_0^{\frac{\pi}{2}}\sin^{2n}tdt=\frac{\pi\binom{2n}{n}}{2^{2n+1}}$, we obtain

$$\displaystyle -2\sum_{n=0}^{\infty}\frac{\zeta(2n)}{\pi^{2n}}\frac{\pi^{2n+2r+1}}{2^{2n+2r+1}(2n+2r)(2n+2r+1)}=\frac{(2r)!\pi}{4^{r+1}}\sum_{n=1}^{\infty}\frac{1}{n^3}\sum_{n_{1}<n_{2}<\ldots n_{r-1}<n}\frac{1}{n_{1}^2n_{2}^2\ldots n_{r-1}^2}$$

or equivalently,

$$\displaystyle -4\sum_{n=0}^{\infty}\frac{\zeta(2n)}{(2n+2r)(2n+2r+1)2^{2n}}\cdot\frac{\pi^{2r}}{(2r)!}=\zeta(\underbrace{2, 2, \ldots, 2}_{\text{$r-1$}}, 3).$$

Now, using the fact that $\displaystyle\zeta(\underbrace{2, 2, \ldots, 2}_{\text{$r$}})=\frac{\pi^{2r}}{(2r+1)!}$ (Proposition 2.3), we obtain the conclusion of the theorem.$\qed{}$
\vspace{0.3cm}

\textbf{Remark.} We believe that for $s>0$ integer, the following equality is true.

\vspace{0.3cm}

\textbf{Conjecture.} \textit{We have 
\bigskip
$$\displaystyle H(r, s)=\frac{-4\pi^{2r+2s+2}}{(2r+2)!}\sum_{n=0}^{\infty}\frac{\zeta(2n)}{(2n+2r+2)(2n+2r+3)\ldots(2n+2r+2s+3)2^{2n}}.$$}
\vspace{0.3cm}

\bigskip

Now, we state the following
\vspace{0.3cm}

\bt[Orr, \cite{Orr}]
For $p\in\mathbb{N}$, and $|z|<1$,

$$\displaystyle\int_0^{\pi z}x^p\cot xdx=(\pi z)^p\sum_{k=0}^{p}\frac{p!(-1)^{[\frac{k+3}{2}]}}{(p-k)!(2\pi z)^{k}}\operatorname{Cl}_{k+1}(2\pi z)+\delta_{\left[\frac{p}{2}\right], \frac{p}{2}}\frac{p!(-1)^{\frac{p}{2}}}{2^p}\zeta(p+1),$$
\bigskip
where $\operatorname{Cl}_{k+1}$ is the $(k+1)$-Clausen function. 
\et
For more basic properties of the Clausen functions, see appendix.

\bigskip

The missing ingredient for our main purpose is given by
\bigskip

\begin{lm}
We have the following equality

$$\displaystyle -2\sum_{n=0}^{\infty}\frac{\zeta(2n)}{(2n+p)2^{2n}}=\log 2+\sum_{k=1}^{[p/2]}\frac{p!(-1)^k(4^k-1)\zeta(2k+1)}{(p-2k)!(2\pi)^{2k}}+\delta_{[p/2], p/2}\frac{p!(-1)^{p/2}\zeta(p+1)}{\pi^p}.$$
\end{lm}
\bigskip

\bigskip

In this moment, we have all the ingredients to prove Zagier's result in the special case ($s=0$),
\bigskip

\textit{Proof of Theorem 2.2}. From Theorem 2.3 we only need to relate the rational zeta series 

$\displaystyle\sum_{n=0}^{\infty}\frac{\zeta(2n)}{(2n+2r+2)(2n+2r+3)2^{2n}}$ to a $\mathbb{Q}$-linear combinations of powers of $\pi$ and odd zeta values. 
\bigskip

Indeed, by applying Lemma 2.6 (for $p=2r+2$ and $p=2r+3$), we have

$$\displaystyle \displaystyle -2\sum_{n=0}^{\infty}\frac{\zeta(2n)}{(2n+2r+2)(2n+2r+3)2^{2n}}=-2\left(\sum_{n=0}^{\infty}\frac{\zeta(2n)}{(2n+2r+2)2^{2n}}-\sum_{n=0}^{\infty}\frac{\zeta(2n)}{(2n+2r+3)2^{2n}}\right)$$

$$\displaystyle=\sum_{k=1}^{r+1}\frac{(2r+2)!(-1)^k(4^k-1)\zeta(2k+1)}{(2r+2-2k)!(2\pi)^{2k}}+\delta_{r+1,r+1}\frac{(2r+2)!(-1)^{r+1}\zeta(2r+3)}{\pi^{2r+2}}$$

$$\displaystyle -\sum_{k=1}^{r+1}\frac{(2r+3)!(-1)^k(4^k-1)\zeta(2k+1)}{(2(r+1-k)+1)!(2\pi)^{2k}}$$

$$\displaystyle =\sum_{k=1}^{r+1}\left[\frac{(2r+2)!(-1)^k(4^k-1)}{(2r+2-2k)!(2\pi)^{2k}}-\frac{(2r+3)!(-1)^k(4^k-1))}{(2(r+1-k)+1)!(2\pi)^{2k}}+\frac{(2r+2)!(-1)^k}{\pi^{2k}}\delta_{r+1, k}\right]\zeta(2k+1)$$

$$\displaystyle =\sum_{k=1}^{r+1}\left[\frac{(-1)^k(4^k-1)}{(2\pi)^{2k}}(2k)!\left(\binom{2r+2}{2k}-\binom{2r+3}{2k}\right)+\delta_{r+1, k}\right]\zeta(2k+1)$$

$$\displaystyle =\sum_{k=1}^{r+1}\frac{(-1)^k(2k)!}{\pi^{2k}}\left[\left(1-\frac{1}{4^k}\right)\left(-\binom{2r+2}{2k-1}\right)+\delta_{r+1, k}\right]\zeta(2k+1)$$

$$\displaystyle =\sum_{k=1}^{r+1}\frac{(-1)^k(2k)!}{\pi^{2k}}\binom{2r+2}{2k-1}\left[\frac{\delta_{r+1, k}}{\binom{2r+2}{2k-1}}-\left(1-\frac{1}{4^k}\right)\right]\zeta(2k+1)$$

$$\displaystyle =\sum_{k=1}^{r+1}\frac{(-1)^k(2k)!}{\pi^{2k}}\cdot\frac{(2r+2)!}{(2k-1)!(2r+2-2k+1)!}\left[\frac{\delta_{r+1, k}}{\binom{2r+2}{2k-1}}-\left(1-\frac{1}{4^k}\right)\right]\zeta(2k+1)$$

$$\displaystyle =\sum_{k=1}^{r+1}\frac{(-1)^k2k}{\pi^{2k}}\cdot\frac{(2r+2)!}{(2(r+1-k)+1)!}\left[\frac{\delta_{r+1, k}}{\binom{2r+2}{2k-1}}-\left(1-\frac{1}{4^k}\right)\right]\zeta(2k+1)$$

$$\displaystyle =\sum_{k=1}^{r+1}\frac{(-1)^k(2k)(2r+2)!\pi^{2(r+1-k)}}{\pi^{2k}\cdot\pi^{2(r+1-k)}\cdot (2(r+1-k)+1)!}\left[\frac{\delta_{r+1, k}}{\binom{2r+2}{2k-1}}-\left(1-\frac{1}{4^k}\right)\right]\zeta(2k+1)$$

$$\displaystyle =\frac{1}{(2r+3)}\sum_{k=1}^{r+1}\frac{(2r+3)!}{\pi^{2r+2}}\cdot\frac{\pi^{2(r+1-k)}}{(2(r+1-k)+1)!}(-1)^k\left[\frac{2k\delta_{r+1, k}}{\binom{2r+2}{2k-1}}-\left(1-\frac{1}{4^k}\right)2k\right]\zeta(2k+1)$$

$$\displaystyle =\frac{1}{(2r+3)\zeta(\underbrace{2, 2, \ldots, 2}_{\text{$r+1$}})}\sum_{k=1}^{r+1}(-1)^k\zeta(\underbrace{2, 2, \ldots, 2}_{\text{$r+1-k$}})c_{r, 0}^k\zeta(2k+1),$$

where obviously $\frac{2k\delta_{r+1,k}}{\binom{2r+2}{2k-1}}=\binom{2k}{2r+2}$ for $1\leq k \leq r+1$. 
\bigskip

By Theorem 1.3, we have

$$\displaystyle H(r, 0):=\zeta(\underbrace{2, 2, \ldots, 2}_{\text{$r$}}, 3)=2(2r+3)\cdot\left(-2\sum_{n=0}^{\infty}\frac{\zeta(2n)}{(2n+2r+2)(2n+2r+3)2^{2n}}\right)\cdot$$
$$\displaystyle\cdot\zeta(\underbrace{2, 2, \ldots, 2}_{\text{$r+1$}}),$$

which will give us

$$\displaystyle H(r, 0)=2(2r+3)\cdot\frac{1}{(2r+3)\zeta(\underbrace{2, 2, \ldots, 2}_{\text{$r+1$}})}\sum_{k=1}^{r+1}(-1)^k\zeta(\underbrace{2, 2, \ldots, 2}_{\text{$r+1-k$}})c_{r, 0}^k\zeta(2k+1)\cdot\zeta(\underbrace{2, 2, \ldots, 2}_{\text{$r+1$}}),$$

which finally reduces to

$$\displaystyle H(r, 0)=2\sum_{k=1}^{r+1}(-1)^kc_{r, 0}^k\zeta(\underbrace{2, 2, \ldots, 2}_{\text{$r+1-k$}})\zeta(2k+1),$$

which is exactly what we wanted to prove.$\qed{}$
\vspace{0.3cm}

\section{A Zagier type formula for multiple $t$-values}
\bigskip

In a similar fashion with the multiple zeta values, we also define the multiple Hurwitz zeta values and multiple $t$-values,

$$\displaystyle\zeta(k_{1}, k_{2}, \ldots, k_{r}; a_{1}, a_{2}, \ldots, a_{r})=\displaystyle\sum_{1\leq n_{1}<n_{2}<\ldots <n_{r}}\frac{1}{(n_{1}+a_{1})^{k_{1}}(n_{2}+a_{2})^{k_{2}}\ldots (n_{r}+a_{r})^{k_{r}}} $$
\bigskip
and

$$\displaystyle t(k_{1}, k_{2}, \ldots, k_{r})=2^{-(k_{1}+k_{2}+\ldots +k_{r})}\zeta(k_{1}, k_{2}, \ldots, k_{r}; -\frac{1}{2}, -\frac{1}{2}, \ldots, -\frac{1}{2})=$$

$$\displaystyle=\sum_{1\leq n_{1}<n_{2}<\ldots<n_{r}}\frac{1}{(2n_{1}-1)^{k_{1}}(2n_{2}-1)^{k_{2}}\ldots (2n_{r}-1)^{k_{r}}}.$$

\bigskip

Unfortunately, there are not many identities known in the literature about multiple Hurwitz zeta values and multiple $t$-values. For more details see \cite{Hoffman2, Murty-Sinha, Shen-Jia, Zhao1} among others. 

Also, let us recall the Gauss hypergeometric function which is defined for $|x|<1$ by the power series

$$\displaystyle {}_{2} F_{1}(a, b; c; x)=\sum_{n=0}^{\infty}\frac{(a)_{n}(b)_{n}}{(c)_{n}}\cdot\frac{x^n}{n!},$$
\vspace{0.3cm}
where $(q)_{n}=q(q+1)\cdot\ldots\cdot (q+n-1)$ is the Pochhammer symbol. As we have already seen in the previous section, the odd integer powers from the Taylor series of $\arcsin$ can be obtained by comparing the coefficients like powers of $\lambda$ in the formulas:

$$\displaystyle\sin(\lambda\arcsin(x))=\lambda x\cdot {}_{2} F_{1}\left(\frac{1+\lambda}{2}, \frac{1-\lambda}{2}; \frac{1}{2}; x^2\right).$$
\bigskip

This implies the following Taylor series for the odd integer powers of arcsin function,

\begin{equation}
   \displaystyle  \displaystyle\frac{\arcsin^{2r+1}(x)}{(2r+1)!}=\sum_{n=0}^{\infty}\frac{\binom{2n}{n}}{(2n+1)4^n}\cdot x^{2n+1}\cdot\sum_{0\leq n_{1}<n_{2}<\ldots<n_{r}<n}\frac{1}{\prod_{i=1}^r(2n_{i}+1)^2}.
\end{equation}

\vspace{0.3cm} 

Again, we begin with the much simpler formula which is given by

\begin{pr} [\cite{Murty-Sinha, Zhao1}] We have the following evaluation

$$t(\underbrace{2, 2, \ldots, 2}_{\text{$r$}})=\frac{\pi^{2r}}{2^{2r}(2r)!}.$$
\end{pr}
\vspace{0.3cm}

\textit{Proof.} Similarly, substitute $x=\sin t$ in formula $(6)$ and we have

$$\displaystyle\frac{t^{2r+1}}{(2r+1)!}=\sum_{n=0}^{\infty}\frac{\binom{2n}{n}}{(2n+1)4^n}\sin^{2n+1}t\sum_{0\leq n_{1}<n_{2}<\ldots<n_{r}<n}\frac{1}{\prod_{i=1}^r(2n_{i}+1)^2}.$$
\vspace{0.3cm}

Again, integrating from $0$ to $\frac{\pi}{2}$ and using Wallis' formula in the form, $\displaystyle\int_0^{\frac{\pi}{2}}\sin^{2n+1}tdt=\frac{1}{2n+1}\cdot\frac{2^{2n}}{\binom{2n}{n}}, n\geq 0$, we have

$$\displaystyle\frac{\pi^{2r+2}}{2^{2r+2}(2r+1)!(2r+2)}=\sum_{n=0}^{\infty}\frac{1}{(2n+1)^2}\sum_{0\leq n_{1}<n_{2}<\ldots<n_{r}<n}\frac{1}{\prod_{i=1}^r(2n_{i}+1)^2}=t(\underbrace{2, 2, \ldots, 2}_{\text{$r+1$}}),$$
\vspace{0.3cm}
and the conclusion follows immediately.$\qed{}$
\bigskip

Next, we evaluate $t(\underbrace{2, 2, \ldots, 2}_{\text{$r$}}, 3)$ in terms of rational zeta series involving $\zeta(2n)$. This formula will be similar with the one produced in Theorem 2.4. 

\bigskip

\bt We have
$$\displaystyle T(r):=t(\underbrace{2, 2, \ldots, 2}_{\text{$r$}}, 3)=-4(r+1)\sum_{n=0}^{\infty}\frac{\zeta(2n)}{(2n+2r+1)(2n+2r+2)2^{2n}}\cdot t(\underbrace{2, 2, \ldots, 2}_{\text{$r+1$}}).$$
\et

\vspace{0.3cm}

\textit{Proof.} Dividing by $x$ and integrating from $0$ to $\sin t$ in formula $(6)$, we have

$$\displaystyle\int_0^{\sin t}\frac{\arcsin^{2r+1}(x)}{x}dx=(2r+1)!\sum_{n=0}^{\infty}\frac{\binom{2n}{n}}{(2n+1)^24^n}\sin^{2n+1}t\sum_{0\leq n_{1}<n_{2}<\ldots<n_{r}<n}\frac{1}{\prod_{i=1}^r(2n_{i}+1)^2}.$$

\vspace{0.3cm}

By the substitution $x=\sin u$ in the integral, we have
\vspace{0.3cm}

$$\displaystyle\int_0^tu^{2r+1}\cot udu=(2r+1)!\sum_{n=0}^{\infty}\sum_{n=0}^{\infty}\frac{\binom{2n}{n}}{(2n+1)^24^n}\sum_{0\leq n_{1}<n_{2}<\ldots<n_{r}<n}\frac{1}{\prod_{i=1}^r(2n_{i}+1)^2}. $$

\vspace{0.3cm}

On the other hand, using $\displaystyle u\cot u=-2\sum_{n=0}^{\infty}\frac{\zeta(2n)}{\pi^{2n}}u^{2n}, |u|<\pi$, we derive

$$\displaystyle\int_0^tu^{2r+1}\cot udu=-2\sum_{n=0}^{\infty}\frac{\zeta(2n)}{\pi^{2n}}\cdot\frac{t^{2n+2r+1}}{2n+2r+1}.$$
\vspace{0.3cm}

Therefore, we obtain
\bigskip

$$\displaystyle -2\sum_{n=0}^{\infty}\frac{\zeta(2n)}{\pi^{2n}}\cdot\frac{t^{2n+2r+1}}{2n+2r+1}=(2r+1)!\sum_{n=0}^{\infty}\frac{\binom{2n}{n}}{(2n+1)^24^n}\sin^{2n+1}t\sum_{0\leq n_{1}<n_{2}<\ldots<n_{r}<n}\frac{1}{\prod_{i=1}^r(2n_{i}+1)^2}.$$

\vspace{0.3cm}

Finally, integrating from $0$ to $\frac{\pi}{2}$ and using Wallis' integral formula in the form, $\displaystyle\int_0^{\frac{\pi}{2}}\sin^{2n+1}tdt=\frac{1}{2n+1}\cdot\frac{2^{2n}}{\binom{2n}{n}}, n\geq 0$, we have
\bigskip

$$\displaystyle\frac{\pi^{2r+2}}{(2r+1)!2^{2r+2}}\left(-2\sum_{n=0}^{\infty}\frac{\zeta(2n)}{(2n+2r+1)(2n+2r+2)2^{2n}}\right)=$$
$$\displaystyle =\sum_{n=0}^{\infty}\frac{1}{(2n+1)^3}\sum_{0\leq n_{1}<n_{2}<\ldots<n_{r}<n}\frac{1}{\prod_{i=1}^r(2n_{i}+1)^2}.$$

\vspace{0.3cm}

Now, using the fact that $t(\underbrace{2, 2, \ldots, 2}_{\text{$r$}})=\frac{\pi^{2r}}{2^{2r}(2r)!}$ (Proposition 3.1) we obtain the desired result.$\qed{}$

\bigskip

\textbf{Remark.} If we define 

$$\displaystyle T(r,s)=t(\{2\}^r,3,\{2\}^s)=t(\underbrace{2, 2, \ldots, 2}_{\text{$r$}}, 3, \underbrace{2, 2, \ldots, 2}_{\text{$s$}}),$$
we believe that the following holds true:
\bigskip

\textbf{Conjecture.} \textit{$$\displaystyle T(r, s)=\frac{-2}{(2r+1)!}\left(\frac{\pi}{2}\right)^{2r+2s+2}\sum_{n=0}^{\infty}\frac{\zeta(2n)}{(2n+2r+1)(2n+2r+2)\ldots(2n+2r+2s+2)2^{2n}}.$$}
\bigskip

In what follows, we derive a Zagier-type formula for the multiple $t$-values in the spirit of Theorem 2.2.
\vspace{0.3cm}

\bt We have

$$\displaystyle T(r)=\sum_{k=1}^{r+1}(-1)^{k+1}d_{r, 0}^k\frac{1}{2^{2k}}\displaystyle\zeta(2k+1)\cdot t(\underbrace{2, 2, \ldots, 2}_{\text{$r+1-k$}}),$$
where $\displaystyle d_{r, 0}^k=\binom{2k}{2r+1}+\left(1-\frac{1}{2^{2k}}\right)2k$, and $r\geq 0$ integer. 

\et

\bigskip

\textit{Proof.} Our intention is to use Lemma 2.6. Indeed, put $p=2r+1$ and $p=2r+2$ and by subtracting, this gives us

$$\displaystyle -2\sum_{n=0}^{\infty}\frac{\zeta(2n)}{(2n+2r+1)(2n+2r+2)2^{2n}}=-2\left(\sum_{n=0}^{\infty}\frac{\zeta(2n)}{(2n+2r+1)2^{2n}}-\sum_{n=0}^{\infty}\frac{\zeta(2n)}{(2n+2r+2)2^{2n}}\right)$$

$$\displaystyle =2\sum_{n=0}^{\infty}\frac{\zeta(2n)}{(2n+2r+2)2^{2n}}-2\sum_{n=0}^{\infty}\frac{\zeta(2n)}{(2n+2r+1)2^{2n}}.$$

$$\displaystyle =\sum_{k=1}^{r+1}\frac{(2r+1)!(-1)^k(4^k-1)\zeta(2k+1)}{(2r+1-2k)!(2\pi)^{2k}}-\sum_{k=1}^{r+1}\frac{(2r+2)!(-1)^k(4^k-1)\zeta(2k+1)}{(2r+2-2k)!(2\pi)^{2k}}$$

$$\displaystyle-\delta_{r+1, r+1}\frac{(2r+2)!(-1)^{r+1}\zeta(2r+3)}{\pi^{2r+2}}$$

$$\displaystyle =\sum_{k=1}^{r+1}\left[\frac{(2r+1)!(-1)^k(4^k-1)}{(2r+1-2k)!(2\pi)^{2k}}-\frac{(2r+2)!(-1)^k(4^k-1)}{(2r+2-2k)!(2\pi)^{2k}}-\delta_{r+1, k}\cdot\frac{(2r+2)!(-1)^k}{\pi^{2k}}\right]\zeta(2k+1)$$

$$\displaystyle =\sum_{k=1}^{r+1}\left[\frac{(-1)^k(4^k-1)}{(2\pi)^{2k}}\left(\frac{(2r+1)!}{(2r+1-2k)!}-\frac{(2r+2)!}{(2r+2-2k)!}\right)-\delta_{r+1, k}\cdot\frac{(2r+2)!(-1)^k}{\pi^{2k}}\right]\zeta(2k+1)$$

$$\displaystyle =\sum_{k=1}^{r+1}\left[\frac{(-1)^k(4^k-1)}{(2\pi)^{2k}}(2k)!\left(\binom{2r+1}{2k}-\binom{2r+2}{2k}\right)-\delta_{r+1, k}\frac{(2r+2)!(-1)^k}{\pi^{2k}}\right]\zeta(2k+1)$$

$$\displaystyle =\sum_{k=1}^{r+1}\frac{(-1)^k(2k)!}{\pi^{2k}}\left[\left(1-\frac{1}{4^k}\right)\left(-\binom{2r+1}{2k-1}\right)-\delta_{r+1, k}\right]\zeta(2k+1)$$

$$\displaystyle =\sum_{k=1}^{r+1}\frac{(-1)^{k+1}(2k)!}{\pi^{2k}}\binom{2r+1}{2k-1}\left[\frac{\delta_{r+1, k}}{\binom{2r+1}{2k-1}}+\left(1-\frac{1}{4^k}\right)\right]\zeta(2k+1)$$

$$\displaystyle =\sum_{k=1}^{r+1}\frac{(-1)^{k+1}(2k)}{\pi^{2k}}\cdot\frac{(2r+1)!}{(2r+2-2k)!}\left[\frac{\delta_{r+1, k}}{\binom{2r+1}{2k-1}}+\left(1-\frac{1}{4^k}\right)\right]\zeta(2k+1)$$

$$\displaystyle =\frac{1}{2r+2}\sum_{k=1}^{r+1}(-1)^{k+1}\frac{(2r+2)!2^{2r+2}}{\pi^{2r+2}}\cdot\frac{\pi^{2(r+1-k)}}{(2(r+1-k))!2^{2(r+1-k)}}\frac{1}{2^{2k}}$$

$$\displaystyle \cdot\left[\frac{(2k)\delta_{r+1, k}}{\binom{2r+1}{2k-1}}+\left(1-\frac{1}{4^k}\right)2k\right]\zeta(2k+1)$$

$$\displaystyle =\frac{1}{(2r+2)t(\underbrace{2, 2, \ldots, 2}_{\text{$r+1$}})}\sum_{k=1}^{r+1}(-1)^{k+1}\left[\binom{2k}{2r+1}+\left(1-\frac{1}{4^k}\right)2k\right]\frac{1}{2^{2k}}\zeta(2k+1)t(\underbrace{2, 2, \ldots, 2}_{\text{$r+1-k$}}).$$

On the other hand, we know by Theorem 3.2 that

$$\displaystyle T(r)=2(r+1)\left(-2\sum_{n=0}^{\infty}\frac{\zeta(2n)}{(2n+2r+1)(2n+2r+2)2^{2n}}\right)\cdot t(\underbrace{2, 2, \ldots, 2}_{\text{$r+1$}}),$$
and finally, we obtain

$$\displaystyle T(r)=\sum_{k=1}^{r+1}(-1)^{k+1}\left[\binom{2k}{2r+1}+\left(1-\frac{1}{4^k}\right)2k\right]\frac{1}{2^{2k}}\zeta(2k+1)t(\underbrace{2, 2, \ldots, 2}_{\text{$r+1-k$}}).$$

$\qed{}$

\section{Appendix} 
\bigskip
The Clausen function (integral) is defined by

$$\displaystyle\operatorname{Cl}_{2}(\theta)=-\int_0^{\theta}\log\left(2\sin\frac{t}{2}\right)dt=\sum_{k=1}^{\infty}\frac{\sin(k\theta)}{k^2},$$

and its Taylor series expansion is given by
\bigskip

$$\displaystyle\frac{\operatorname{Cl}_{2}(\theta)}{\theta}=1-\log|\theta|+\sum_{n=1}^{\infty}\frac{\zeta(2n)}{n(2n+1)}\left(\frac{\theta}{2\pi}\right)^{2n}, |\theta|<2\pi.$$
\bigskip

The higher order Clausen functions are

$$\displaystyle\operatorname{Cl}_{2m}=\sum_{k=1}^{\infty}\frac{\sin(k\theta)}{k^{2m+1}}, \operatorname{Cl}_{2m+1}(\theta)=\sum_{k=1}^{\infty}\frac{\cos(k\theta)}{k^{2m+1}}.$$
\bigskip

Using the properties of the Riemann zeta function, we have the following particular values:

$$\displaystyle\operatorname{Cl}_{2m}(\pi)=0, \operatorname{Cl}_{2m+1}(\pi)=-\frac{(4^m-1)\zeta(2m+1))}{4^m} ,$$

and

$$\displaystyle\operatorname{Cl}_{2m}\left(\frac{\pi}{2}\right)=\beta(2m), \operatorname{Cl}_{2m+1}\left(\frac{\pi}{2}\right)=-\frac{(4^m-1)\zeta(2m+1)}{2^{4m+1}},$$

where $\beta(s)=\sum_{n=0}^{\infty}\frac{(-1)^n}{(2n+1)^s}, \operatorname{Res}>0$ is the Dirichlet beta function.

Moreover,

$$\displaystyle\frac{d}{d\theta}\operatorname{Cl}_{2m}(\theta)=\operatorname{Cl}_{2m-1}(\theta), \frac{d}{d\theta}\operatorname{Cl}_{2m+1}(\theta)=-\operatorname{Cl}_{2m}(\theta),$$

and

$$\displaystyle\int_0^{\theta}\operatorname{Cl}_{2m}(x)dx=\zeta(2m+1)-\operatorname{Cl}_{2m+1}(\theta), \int_0^{\theta}\operatorname{Cl}_{2m-1}(x)dx=\operatorname{Cl}_{2m}(\theta).$$
\bigskip

\bigskip

\textit{Proof of the Lemma 2.5} (see also \cite{Orr}). Theorem 2.4 for $z=\frac{1}{2}$ gives us

$$\displaystyle\int_0^{\frac{\pi}{2}}x^p\cot xdx=\left(\frac{\pi}{2}\right)^p\left(\log 2+\sum_{k=1}^{[p/2]}\frac{p!(-1)^k(4^k-1)}{(p-2k)!(2\pi)^{2k}}\zeta(2k+1)\right)+\delta_{\left[\frac{p}{2}\right], \frac{p}{2}}\frac{p!(-1)^{\frac{p}{2}}\zeta(p+1)}{2^p}.$$

On othe other hand, since $\displaystyle\cot x=-2\sum_{n=0}^{\infty}\frac{\zeta(2n)}{\pi^{2n}}\cdot x^{2n-1}, |x|<\pi$, by integration and Fubini's theorem, we obtain

$$\displaystyle\int_0^{\frac{\pi}{2}}x^p\cot xdx=\int_0^{\frac{\pi}{2}}x^p\left(-2\sum_{n=0}^{\infty}\frac{\zeta(2n)}{\pi^{2n}}\cdot x^{2n-1}\right)dx=-2\sum_{n=0}^{\infty}\frac{\zeta(2n)}{\pi^{2n}}\int_0^{\frac{\pi}{2}}x^{2n+p-1}dx$$

or equivalently,

$$\displaystyle\int_0^{\frac{\pi}{2}}x^p\cot xdx=-2\left(\frac{\pi}{2}\right)^p\sum_{n=0}^{\infty}\frac{\zeta(2n)}{(2n+p)2^{2n}},$$

and the lemma follows immediately.$\square$

\bigskip

\begin{flushright}
\begin{minipage}{148mm}\sc\footnotesize
Texas Tech University, Department of Mathematics \& Statistics, 1108 Memorial Circle, Lubbock, TX 79409, USA\\
{\it E--mail address}: {\tt Cezar.Lupu@ttu.edu, lupucezar@gmail.com} \vspace*{3mm}
\end{minipage}
\end{flushright}


\bigskip

\bigskip


\begin{thebibliography}{199}
\markboth{Bibliography}{Bibliography}


\bibitem{Andrews-Askey-Roy}
G. E. Andrews, R. Askey, R. Roy, {\em Special Functions}, Cambridge University Press, 2000.

\bibitem{Ayoub}
R. Ayoub, Euler a d the zeta function, {\em Amer.Math. Monthly} \textbf{81} (1974), 1067--1086.

\bibitem{Borwein-Broadhurst-Crandall}
J. M. Borwein, D. Braodhurst, R. E. Crandall, Computational strategies for the Riemann zeta function, {\em J. Comp. Appl. Math.} \textbf{121} (2000), 247--296. 

\bibitem{Borwein-Chamberland}
J. Borwein, M. Chamberland, Integer powers of arcsin, {\em Int. J. Math. Math. Sci.} \textbf{2007} (2007), Article ID 19381, 10pp.

\bibitem{Brown}
F. Brown, Mixed Tate motives over $\mathbb{Z}$, {\em Annals of Math.} \textbf{175} (2012), 949--976.

\bibitem{Brown-Schnetz}
F. Brown, O. Schnetz, Single-valued multiple polylogarithms and a proof of the zig-zag conjecture, {\em J. Number Theory} \textbf{148} (2015), 478--506.

\bibitem{Burgos-Fresan}
J. I. Burgos Gil, J. Fresan, Multiple zeta values: from numbers to motives, {\em Clay Mathematics Proceedings}, to appear.

\bibitem{Charlton}
S. Charlton, $\zeta(\{\{2\}^m, 1, \{2\}^m, 3\}^n, \{2\}^m)/\pi^{4n+2m(2n+1)}$ is irrational, {\em J. Number Theory} \textbf{148} (2015), 463--477.

\bibitem{Chu-Zheng}
W. Chu, D. Zheng, Infinite series with harmonic numbers and central binomial coefficients, {\em Int. J. Number Theory} \textbf{5} (2009), 429--448.

\bibitem{Chung}
C.L. Chung, On the sum relation of multiple Hurwitz zeta functions, {\em Quest. Math.}, to appear.

\bibitem{Clausen}
T. Clausen,  Uber die function $\sin\phi+\frac{1}{2^2}\sin 2\phi+\frac{1}{3^2}\sin 3\phi+$etc., {\em J. Reine Angew. Math.} \textbf{8} (1832), 298–-300.

\bibitem{Hoffman}
M. E. Hoffman, Multiple harmonic sums, {\em Pacific J. Math.} \textbf{152} (1992), 275--290.

\bibitem{Hoffman1}
M. Hoffman, The algebra of multiple harmonic series, {\em J. Algebra} \textbf{194} (1997), 477–495.

\bibitem{Hoffman2}
M. Hoffman, An odd variant of multiple zeta values, {\em Commun. Number Theory Phys.} \textbf{13} (2019), 529–567.

\bibitem{Leshchiner}
D. Leshchiner, Some new identities for $\zeta(k)$, {\em J. Number Theory} \textbf{13} (1981), 355--362.

\bibitem{Li}
Z-H. Li, Another proof of Zagier's evaluation formula of the multiple zeta values $\zeta(2, \ldots, 3, 2, \ldots, 2)$, {\em Math. Res. Lett.} \textbf{20} (2013), 947--950.

\bibitem{Lupu}
C. Lupu, {\em Analytic Aspects of the Riemann Zeta and Multiple Zeta Values}, PhD Thesis, available at $\mathtt{http://d-scholarship.pitt.edu/35330/}$, University of Pittsburgh, 2018.

\bibitem{Murty-Sinha}
M. Ram Murty, K. Sinha, Multiple Hurwitz zeta functions, {\em Proc. Symposia Pure Math.} \textbf{75} (2006), 135--156.

\bibitem{Orr}
D. Orr, Generalized rational zeta series for $\zeta(2n)$ and $\zeta(2n+1)$, {\em Integral Transforms Spec. Func.} \textbf{28} (2017), 966--987. 

\bibitem{Rainville}
E. D. Rainville, {\em Special Functions}, Chelsea Publishing House, 1971.

\bibitem{Shen-Jia}
Z. Shen, L. Jia, Some identities for multiple Hurwitz zeta values, {\em J. Number Theory} \textbf{179} (2017), 256--267.

\bibitem{Sun}
Z. W. Sun, New series for some special values of $L$-functions, {\em Nanjiang Univ. J. Math. Biquarterly} \textbf{32} (2015), 189--218.

\bibitem{ZagierMZV}
D. Zagier, Values of zeta functions and their applications, {\em Proceedings of the First European Congress of Mathematicians}, 1992. 

\bibitem{Zagier}
D. Zagier, Multiple zeta values, {\em Unpublished manuscript}, Bonn, 1995.

\bibitem{Zagier1}
D. Zagier, Evaluation of the multiple zeta values $\zeta(2, \ldots, 2, 3, 2, \ldots, 2)$, {\em Ann. Math.} \textbf{175} (2012), 977--1000. 

\bibitem{Zhao1}
J. Zhao, Sum formula of multiple Hurwitz-zeta values, {\em Forum Math.} \textbf{27} (2015), 929--936.

\bibitem{Zhao}
J. Zhao, {\em Multiple Zeta Functions, Multiple polylogarithms, and their Special Values}, World Scientific Publishing Co., 2016.




\end{thebibliography}
\end{document}